\input amstex
\magnification=\magstep0
\documentstyle{amsppt}
\pagewidth{6.0in}
\input amstex
\topmatter
\title
A note on Fibonacci-type polynomials
\endtitle
\author
Tewodros Amdeberhan
\endauthor
\affil
Department of Mathematics\\
Tulane University\\
tamdeber\@tulane.edu
\endaffil
\abstract 
We opt to study the convergence of maximal real roots of certain Fibonacci-type polynomials given by $G_n=x^kG_{n-1}+G_{n-2}$. The special cases $k=1$ and $k=2$ are found in [4] and [7], respectively.
\endabstract
\endtopmatter
\def\({\left(}
\def\){\right)}

\document
\noindent
In the sequel, $\Bbb{P}$ denotes the set of positive integers. 
The Fibonacci polynomials [2] are defined recursively by 
$F_0(x)=1, F_1(x)=x$ and
$$F_n(x)=xF_{n-1}(x)+F_{n-2}(x),\qquad n\geq 2.$$
\bf Fact 1. \it Let $n\geq 1$. Then the roots of $F_n(x)$ are given by 
$$x_k=2i\cos\left(\frac{\pi k}{n+1}\right),\qquad 1\leq k\leq n.$$
In particular a Fibonacci polynomial has no positive real roots. \rm
\smallskip
\noindent
\bf Proof. \rm The Fibonacci polynomials are essentially Tchebycheff
polynomials. This is well-known (see, for instance [2]). $\square$
\bigskip
\noindent
Let $k\in\Bbb{P}$ be fixed. Several authors ([3]-[7]) have investigated the so-called \it Fibonacci-type \rm polynomials. In this note, we focus on a particular group of polynomials recursively defined by
$$G_n^{(k)}(x)=
\cases
-1, \qquad &n=0\\x-1, \qquad &n=1\\
x^k G_{n-1}^{(k)}(x) + G_{n-2}^{(k)}(x),\qquad &n\geq 2.\endcases$$
When there is no confusion, suppress the index $k$ to write $G_n$ for $G_n^{(k)}(x)$. We list a few basic properties relevant to our work here.
\smallskip
\noindent
\bf Fact 2. \it For each $k\in\Bbb{P}$, there is a rational generating function for $G_n$; namely,
$$\sum_{n\geq 0}G_n^{(k)}(x)t^n=\frac{(x^k+x-1)t-1}{1-x^kt-t^2}.$$ \rm
\bf Proof. \rm follows from the definition of $G_n$. $\square$
\smallskip
\noindent
\bf Fact 3. \it The following relation holds
$$G_n^{(k)}(x)=\frac{G_{n-1}^{(k)}F_{n-1}(x^k)+(-1)^{n-1}}{F_{n-2}(x^k)}.$$ \rm
\bf Proof. \rm Write the equivalent formulation

$$G_n^{(k)}(x)=\det\pmatrix x-1 & -1 & 0 & 0 &\hdots & 0&0&0\\
            -1 & x^k&-1 & 0 &\hdots & 0&0&0\\
             0 & 1 & x^k& -1&\hdots & 0&0&0\\
\hdotsfor 8\\
             0 & 0 &  0 & 0&\hdots & x^k &-1&0\\
             0 & 0 &  0 & 0&\hdots &1& x^k &-1\\
             0 & 0 &  0 & 0&\hdots &0 &1 & x^k
\endpmatrix,$$ 
then apply Dodgson$^{\prime}$s determinantal formula [1]. $\square$
\bigskip
\noindent
\bf Fact 4. \it 
For a fixed $k$, let $\{g_n^{(k)}\}_{n\in\Bbb{P}}$ be the maximal real roots of $\{G_n^{(k)}(x)\}_{n}$. Then $\{g_{2n}^{(k)}\}_n$ and $\{g_{2n-1}^{(k)}\}_n$ are decreasing and increasing sequences, respectively. \rm
\smallskip
\noindent
\bf Proof. \rm First, each $g_n$ exists since $G_n(0)=1<0$ and $G_n(\infty)=\infty$. Assume $x>0$. Invoking Fact 3 from above, twice, we find that
$$F_{2n-3}(x^k)G_{2n}^{(k)}(x)=F_{2n-1}(x^k)G_{2n-2}^{(k)}(x)+x^k, \qquad
F_{2n-2}(x^k)G_{2n+1}^{(k)}(x)=F_{2n}(x^k)G_{2n-1}^{(k)}(x)-x^k.$$
From these equations and $F_n(x)>0$ (see Fact 1), it is clear that
$G_{2n-2}(x)>0$ implies $G_{2n}(x)>0$; also if $G_{2n-2}(x)=0$ then $G_{2n}(x)>0$. Thus $g_{2n-2}>g_{2n}$. A similar argument shows $g_{2n+1}>g_{2n-1}$. The proof is complete. $\square$
\bigskip
\noindent
Define the quantities
$$\align
\alpha(x)&=\frac{x+\sqrt{x^2+4}}2,\qquad {}\qquad \beta(x)=\frac{x-\sqrt{x^2+4}}2,\\
p(x)&=\frac{(x-1)+\beta(x^k)}{\alpha(x^k)-\beta(x^k)}, \qquad
q(x)=\frac{(x-1)+\alpha(x^k)}{\alpha(x^k)-\beta(x^k)}.\endalign$$
\bf Fact 5. \it For $n\geq 0$ and $k\in\Bbb{P}$, we have the explicit formula
$$G_n^{(k)}(x)=p(x)\alpha^n(x^k)-q(x)\beta^n(x^k).$$
\bf Proof. \rm this is a standard procedure. $\square$

$$ $$

\pagebreak

\noindent
For each $k\in\Bbb{P}$, introduce another set of polynomials
$$H^{(k)}(x)=x^k-x^{k-1}+x-2.$$
\bf Fact 6. \it For each $k\in\Bbb{P}$, the polynomial $H^{(k)}(x)$ has exactly one positive real root $\xi^{(k)}$. And $\xi^{(k)}> 1$. \rm
\smallskip
\noindent
\bf Proof. \rm Since $H^{(k)}(x)=(x-1)(x^{k-1}+1)-1<0$, whenever $0<x\leq 1$, there are no roots in the range $0<x\leq 1$. 
On the other hand, 
$H^{(k)}(1)<0, H^{(k)}(\infty)=\infty$ and the derivative 
$$\frac{d}{dx}H^{(k)}(x)=x^{k-1}(k(x-1)+1)+1> 0 \qquad\text{whenever $x\in\Bbb{P}$},$$
suggest there is only one positive root (necessarily greater than 1). $\square$
\bigskip
\noindent
\bf Fact 7. \it If $k$ is odd (even), then $H^{(k)}(x)$ has no (exactly one) negative real root. \rm
\smallskip
\noindent
\bf Proof. \rm For $k$ odd, $H^{(k)}(-x)=(-x-1)(x^{k-1}+1)-1<0$. For $k$ even, $H^{(k)}(-x)=x^k+x^{k-1}-x-2$ changes sign only once. Apply Descarte$^{\prime}$s Rule. $\square$
\smallskip
\noindent
Now, we state and prove the main result of the present note.
\bigskip
\noindent
\bf Theorem. \it Preserve the notations of Facts 4 and 6. Then, depending on the parity of $n$, the roots $\{g_n^{(k)}\}_n$ converge from above or below so that
$g_n^{(k)}\rightarrow\xi^{(k)}$ as $n\rightarrow\infty$. Note also $\xi^{(k)}\rightarrow 1$ as $k\rightarrow\infty$. \rm
\bigskip
\noindent
\bf Proof. \rm For notational brevity, suppress $k$ and write $g_n$ and $\xi$.
From $G_n(g_n)=0$ and Fact 5 above, we resolve
$$\frac{p(g_n)}{q(g_n)}=\frac{\beta^n(g_n^k)}{\alpha^n(g_n^k)},\qquad\text{or}
\qquad
\frac{2(g_n-1)+g_n^k-\sqrt{g_n^{2k}+4}}{2(g_n-1)+g_n^k+\sqrt{g_n^{2k}+4}}=
(-1)^n\left(1-\frac{2g_n}{g_n^k+\sqrt{g_n^{2k}+4}}\right)^n.\tag1$$
Using Gershgorin$^{\prime}$s Circle theorem, it is easy to see that $1\leq g_n\leq 2$. When combined with Fact 4, the monotonic sequences $\{g_{2n}\}_n$ and $\{g_{2n-1}\}_n$ converge to finite limits, say $\xi_{+}$ and $\xi_{-}$ respectively. 
\smallskip
\noindent
The right-hand side of (1) vanishes in the limit $n\rightarrow\infty$, thus
$$2(\xi-1)+\xi^k-\sqrt{\xi^{2k}+4}=0.$$
Further simplification leads to $H^{(k)}(\xi)=\xi^k-\xi^{k-1}+\xi-2=0$. From Fact 6, such a solution is unique. So, $\xi_{+}=\xi_{-}=\xi$ completes the proof. $\square$

$$ $$

\pagebreak
\bigskip
\noindent
\bf Appendix \rm
\bigskip
\noindent
In this section, we discuss the \it bivariate Fibonacci \rm polynomials,
of the \it first kind \rm (BFP1), defined as
$$g_n(x,y)=xg_{n-1}(x,y)+yg_{n-2}(x,y),\qquad g_0(x,y)=x, \qquad
g_1(x,y)=y.$$
If $x=y=1$ then the resulting sequence is the Fibonacci numbers. 
\bigskip
\noindent
The following is a generating function for the BVP1
$$\sum_{n\geq 0}g_n(x,y)t^n=\frac{x+(y-x^2)t}{1-xt-yt^2}.$$
It is also possible to give an explicit expression
$$g_n(x,y)=\sum_{k\geq 1}
\frac{2n-3k+1}{n-k}\binom{n-k}{k-1}x^{n-2k+1}y^k.$$
This shows clearly that each BFP1 has non-negative coefficients. 
\bigskip
\noindent
The other variant appears often in the literature which we call
\it bivariate Fibonacci \rm polynomials, of the \it second kind \rm 
(BFP2). These are recursively defined as
$$f_n(x,y)=xf_{n-1}(x,y)+yf_{n-2}(x,y),\qquad f_0(x,y)=y, \qquad
f_1(x,y)=x.$$
Obviously $f_n(1,1)$ yields the Fibonacci numbers. We also find the 
ordinary generating function
$$\sum_{n\geq 0}f_n(x,y)t^n=\frac{y+(x-xy)t}{1-xt-yt^2}.$$
One interesting contrast between the two families is the following.
While the roots of $f_n(x,1)$ are all imaginary, the roots of
$g_n(1,y)$ are all real numbers.
\bigskip
\noindent
Using the corresponding generating functions for BVP2 $f_n(x,y)$ and the
classical Fibonacci polynomials $F_n(x)=f_n(x,1)$ proves the below affine 
relation
$$f_n(x,y^2)=xy^{n-1}F_{n-1}(x/y)+y^{n+2}F_{n-2}(x/y).$$
In particular, the \it Jacobsthal-Lucas \rm numbers $J_n=f_n(2,1)$ can be expreesed in terms of 
values of the Fibonacci polynomials, at $1/\sqrt{2}$, namely that
$$J_n=2^{\frac{n-1}2}F_{n-1}\left(\frac1{\sqrt{2}}\right)+
2^{\frac{n}2+1}F_{n-2}\left(\frac1{\sqrt{2}}\right).$$
Since we have
$$\sum_{n\geq 0}F_n(x)t^n=\frac{1}{1-xt-t^2}\qquad\text{and}\qquad
F_n(x)=\sum_{k\geq 0}\binom{n-k}kx^{n-2k},$$
we obtain
$$f_n(x,y^2)=\sum_{k\geq 0}\binom{n-k-1}kx^{n-2k}y^k+
\sum_{k\geq 0}\binom{n-k-2}kx^{n-2k-2}y^{k+2}.$$
In particular, when $x=1$ there holds
$$f_n(1,y)=\sum_{k=0}^{\lfloor{(n+1)/2}\rfloor}
\frac{(n-k-1)!}{k!(n-2k+2)}Q(n,k)y^k$$
where $Q(n,k)=n^3-3(2k-1)n^2+(13k(k-1)+2)n-k(k-1)(9k-4)$. 
\bigskip
\noindent
If we alter the definition of BFP2 and specialize as $h_0=2, h_1=1, 
h_n(x)=h_{n-1}(x)+xh_{n-2}(x)$ then the resulting
sequence of polynomials become intimately linked to the Lucas polynomials $L_n(x)$ as follows
$$L_n(x)=x^nh_n(1/x^2).$$

\enddocument